\documentclass[12pt]{article}
\usepackage{mathrsfs}
\usepackage{amsmath,amsfonts,amssymb,rotating,psfrag}
\usepackage{hyperref}
\usepackage{amsfonts}
\usepackage{latexsym,amscd,epsfig,graphicx,colordvi}
\usepackage{ifpdf}
\usepackage{xcolor}
\usepackage[dvipsnames]{pstricks}

\def\qed{\nopagebreak\hfill{\rule{4pt}{7pt}}}

\newtheorem{theo}{Theorem}[section]

\newtheorem{defi}[theo]{Definition}
\numberwithin{equation}{section}

\newdimen\Squaresize \Squaresize=11pt
\newdimen\Thickness \Thickness=0.7pt
\def\Square#1{\hbox{\vrule width \Thickness
   \vbox to \Squaresize{\hrule height \Thickness\vss
    \hbox to \Squaresize{\hss#1\hss}
   \vss\hrule height\Thickness}
\unskip\vrule width \Thickness} \kern-\Thickness}

\def\Vsquare#1{\vbox{\Square{$#1$}}\kern-\Thickness}

\def\moins{\raise 1pt\hbox{{$\scriptstyle -$}}}
\title{Overpartitions Identities Involving Gaps and Weights}
\author{Jeremy J.F. Guo\raisebox{5pt}{\scriptsize 1},
Doris D. M. Sang\raisebox{5pt}{\scriptsize 2} and Diane Y. H.
Shi\raisebox{5pt}{\scriptsize 3}}
\date{
\vspace{15pt} \raisebox{5pt}{\scriptsize 1\,}Center for Applied Mathematics\\School of
Mathematics\\ Tianjin University\\ Tianjin 300072, P.R.
China\\[1mm] Email: shiyahui@tju.edu.cn\\
\vspace{15pt}\raisebox{5pt}{\scriptsize 2\,}School of Mathematics and Quantitative Economics\\Dongbei
University of Finance and Economics\\ Liaoning 116025, P.R. China\\[1mm] Email:  sangdm@dufe.edu.cn\\
\vspace{15pt}\raisebox{5pt}{\scriptsize 3\,}Center for Applied Mathematics\\ Tianjin University\\School of
Mathematics\\ Tianjin 300072, P.R.
China\\[1mm] Email: shiyahui@tju.edu.cn}

\begin{document}

\allowdisplaybreaks

\maketitle

\vspace{0.3cm} \noindent{\bf Abstract.}

In this paper, we find an identity which connects the overpartition function and the function of Rogers--Ramanujan--Gordon type overpartitions by considering the weights and gaps. This identity can be seen as an analogue of the  weighted identity of Alladi.



\section{Introduction}

In 1997, Alladi \cite{all} began a systematic study of weighted partition identities.
By attaching weights to the gaps between parts of partitions, he obtained interesting new identities connecting the famous partition functions of Euler, Gauss, Lebesgue, Rogers--Ramanujan and others. In this paper we apply the definition of the successive Durfee squares for overpartitions given by Corteel and Mallet \cite{cor07} to give
 an identity connecting the Rogers--Ramanujan-Gordon type overpartitions and the overpartitions without any restrictions, which can be seen as an analogue of the  weighted identity of Alladi.

To illustrate
 the weighted partition identities he obtained, Alladi proposed the following general problem: Given a set $S$ of partitions, let $P_S(n)$ denote the number of partitions $\pi$ of $n$ with $\pi\in S$. Suppose $S\subset T$. The problem is to determine weights $w_S(\pi)\geq 1$ in a natural way such that for all $n$
\begin{equation}\label{pgw}
P_T(n)=\sum_{\sigma(\pi)=n,\ \pi\in S}w_S(\pi).
\end{equation}
(Here and in what follows, $\sigma(\pi)$ is the sum of the parts of $\pi$.)

 Naturally, we want to know whether we could get some new identities, if the sets $S$ and $T$ in above problem are replaced by sets of overpartitions. We focus on the situation that  $S$ denotes the set of the Rogers--Ramanujan-Gordon type overpartitions and $T$ denotes the set of the overpartitions without any restrictions.

Recall that a partition $\lambda$ of a positive integer $n$ is a finite non-increasing sequence of positive integers $(\lambda_1,
\lambda_2,\ldots,\lambda_l)$ such that $\sum_{i=1}^{l}\lambda_i=n$. Let 
$p(n)$ denote the number of partitions of $n$.

An overpartition is a partition for which the first occurrence of a part may be overlined. For example, $(\overline{7},7,6,\overline{5},2,\overline{1})$ is an overpartition of $28$. The number of overpartitions of $n$ is denoted by the function $\overline{p}(n)$.

For a partition $\lambda$, we use the same symbol to denote its Ferrers graph. The conjugate of the Ferrers graph of $\lambda$ is denoted by $\lambda^{*}$.

Given two partitions $\lambda_1$ and $\lambda_2$, by $\lambda_1+\lambda_2$ we mean the partition whose parts are obtained by adding the number of nodes in the corresponding rows of the Ferrers graphs of $\lambda_1$ and $\lambda_2$. If $\lambda_1$ is an overpartition, then the parts of $\lambda_1+\lambda_2$ in the corresponding rows of the Ferrers graph of overlined parts of $\lambda_1$ remain overlined.

Here and in the rest of this paper, we employ the customary q-series notations
\[(a)_\infty=(a;q)_{\infty}=\sum_{i=0}^{\infty}(1-aq^i),\]

\[(a)_n=(a;q)_n=\frac{(a)_{\infty}}{(aq^n)_{\infty}}\]
and
\[(a_1,\ldots,a_k;q)_{\infty}=(a_1;q)_{\infty}\cdots(a_k;q)_{\infty}.\]

By  these notations, we recall the generating functions of the partitions and overpartitions
\begin{equation}
\sum_{n}p(n)q^n=\frac{1}{(q;q)_\infty},
\end{equation}
\begin{equation}
  \sum_{n}\overline{p}(n)q^n=\frac{(-q)_\infty}{(q)_\infty}.
\end{equation}

By a combinatorial proof,
Corteel and Lovejoy \cite{cor04} gave another form of the generating function of overpartitions
\begin{equation}\label{2over}
\sum_{n}\overline{p}(n)q^n=\frac{(-q)_\infty}{(q)_\infty}=\sum_{n=0}^\infty\frac{(-1;q)_nq^{n(n+1)/2}}{(q;q)^2_n},
\end{equation}
which will be employed to demonstrate the relation between the overpartitions without any restrictions and the Rogers-Ramanujan-Gordon type overpartitions.

The theorem of the general Rogers-Ramanujan-Gordon type overpartitions was given by  Chen, Sang and Shi \cite{chen} which can be seen as an analogue of  Gordon's \cite{gor61} combinatorial generalization of the Rogers--Ramanujan identities. The definition of the Rogers--Ramanujan--Gordon type overpartition is subtracted from this theorem which is stated as follows.
 \begin{defi}\label{dki}For $k\geq 2$ and  $k\geq i\geq 1$, we define the Rogers--Ramanujan--Gordon type overpartition $\lambda$ to be the overpartition such that $\lambda_1+\lambda_2+\cdots+\lambda_s$ with part $1$  occurs as a non-overlined part at most $i-1$ times,
and  $\lambda_j-\lambda_{j+k-1}\geq1$ if $\lambda_j$ is overlined and
$\lambda_j-\lambda_{j+k-1}\geq2$ otherwise.
Let $D_{k,i}(n)$ denote the number of such overpartitions $\lambda$ of $n$.\end{defi}

Chen, Sang and Shi \cite{chen} also derived the following generating function form of $D_{k,i}(n)$ as well.
For $k\geq 2$ and  $k\geq i\geq 1$, we have
\begin{align}\nonumber
&\sum_{n=0}^{\infty}D_{k,i}(n)q^n\\&\qquad\label{zki}=\sum_{N_1\geq N_2\geq\cdots\geq N_{k-1}
\geq0}\frac{q^{N_1(N_1+1)/2+N_2^2+\cdots+N_{k-1}^2+N_{i+1}+\cdots+N_{k-1}}(-q)_{N_1-1}(1+q^{N_i})}
{(q)_{N_1-N_2}\cdots(q)_{N_{k-2}-N_{k-1}}(q)_{N_{k-1}}}.
\end{align}
Setting $i=k$ in \eqref{zki}, we get the generating function of $D_{k,k}(n)$
\begin{align}\nonumber
&\sum_{n=0}^{\infty}D_{k,k}(n)q^n\\&\qquad\label{zkk}=\sum_{N_1\geq N_2\geq\cdots\geq N_{k-1}
\geq0}\frac{q^{N_1(N_1+1)/2+N_2^2+\cdots+N_{k-1}^2}(-1)_{N_1}}
{(q)_{N_1-N_2}\cdots(q)_{N_{k-2}-N_{k-1}}(q)_{N_{k-1}}}.
\end{align}

Applying the definition of successive Durfee squares for overpartitions, we could give new combinatorial explanations for \eqref{2over} and \eqref{zkk}. Recall that in 1979, Andrews \cite{and79}
introduced the idea of Durfee dissection of an integer
 partition to interpret  the generalized Rogers-Ramanujan identity combinatorially. To extend this work to overpartitions, Corteel and Mallet in \cite{cor07} introduced the generalized Durfee square for overpartitions.

In Section 2, we shall introduce the definition of successive Durfee squares for overpartitions, and give new combinatorial explanations for \eqref{2over} and \eqref{zkk}. Then in Section 3, we will apply these new explanations to give the definition of $w_S(\pi)$ where $P_T(n)=\overline{p}(n)$ and $P_S(n)=D_{k,k}(n)$ in \eqref{pgw}.


\section{Successive Durfee Squares}

In this section, we shall introduce the definition of successive Durfee squares for overpartitions, and give new combinatorial explanations for \eqref{2over} and \eqref{zkk}.

Firstly, we will introduce the definition of generalized Durfee square.

 \begin{defi}\label{gds}(See \cite{cor07}.)
We say that the generalized Durfee square of an overpartition $\lambda$ has size $N$ if $N$ is the largest integer
  such that the number of overlined parts  plus
the number of non-overlined parts greater than or equal to $N$   is greater than or equal
to $N$. Let $D(\alpha)$ denote the size of generalized Durfee square of $\alpha$.
\end{defi}

For example, let $\alpha=(7,6,6,\overline{5},\overline{3}, 3,2,\overline{1})$. We rewrite it as $\alpha=(\overline{5},\overline{3},\overline{1},7,6,6, 3,2)$. And we could get that $D(\alpha)=6$. Therefore, $\alpha$ has a generalized Durfee square of size $6$. We present its Ferrers graph as follows.

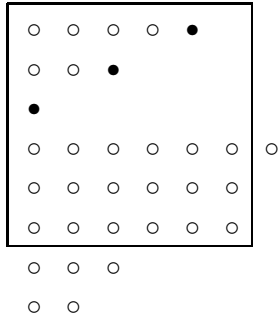
\begin{center}
 \begin{picture}(100,130)
  \put(10,15){\circle{4}}
  \put(25,15){\circle{4}}

  \put(10,30){\circle{4}}
  \put(25,30){\circle{4}}
  \put(40,30){\circle{4}}

  \put(10,45){\circle{4}}
  \put(25,45){\circle{4}}
  \put(40,45){\circle{4}}
\put(55,45){\circle{4}}
  \put(70,45){\circle{4}}
  \put(85,45){\circle{4}}

  \put(10,60){\circle{4}}
  \put(25,60){\circle{4}}
  \put(40,60){\circle{4}}
  \put(55,60){\circle{4}}
  \put(70,60){\circle{4}}
  \put(85,60){\circle{4}}

  \put(10,75){\circle{4}}
  \put(25,75){\circle{4}}
  \put(40,75){\circle{4}}
  \put(55,75){\circle{4}}
  \put(70,75){\circle{4}}
  \put(85,75){\circle{4}}
  \put(100,75){\circle{4}}

  \put(10,90){\circle*{4}}

  \put(10,105){\circle{4}}
  \put(25,105){\circle{4}}
  \put(40,105){\circle*{4}}

  \put(10,120){\circle{4}}
  \put(25,120){\circle{4}}
  \put(40,120){\circle{4}}
\put(55,120){\circle{4}}
  \put(70,120){\circle*{4}}

\put(0,38){\line(0,1){92}}

\put(0,38){\line(1,0){92.5}}

\put(92.5,38){\line(0,1){92}}

\put(0,130){\line(1,0){92.5}}
\end{picture}
\end{center}

Inspired by the Algorithm Z \cite{and84}, Corteel and Mallet in \cite{cor07} also gave the generating function for ovrpartitions with generalized Durfee square of size $N$ where the exponent of $q$ counts the weight and the exponent of $a$ counts  the number of overlined parts as follows
\begin{equation}
  \frac{a^Nq^{\frac{(N+1)N}{2}}(-1/a)_N}{(q)_N(q)_N}.
\end{equation}

Following this idea, we can use the definition of  generalized Durfee square to explain the common factor both in \eqref{2over} and \eqref{zkk}
\begin{equation}
  \sum_{n=0}^{\infty}\frac{(-1;q)_nq^{(n+1)n/2}}
  {(q;q)_n}.
\end{equation}

\begin{theo}
Let $\alpha$ denote the overpartition such that the number of parts of $\alpha$ equals the size of
generalized Durfee square of $\alpha$. And $g(n)$
 denotes the number of overpartitions $\alpha$ of $n$.
Then the generating function of $g(n)$ is as follows
\begin{equation}
  \sum_{n=0}^{\infty}g(n)q^n=\sum_{N=0}^{\infty}\frac{(-1;q)_Nq^{(N+1)N/2}}
  {(q;q)_N}.
\end{equation}

\end{theo}


\noindent Proof. Let $\gamma=(\gamma_1,\gamma_2,\ldots,\gamma_N)$
denote a  partition with $N$ distinct parts and $\delta=(\delta_1,\delta_2,\ldots,\delta_k)$ denote a  partition with distinct
parts and $0\leq\delta_i\leq N-1$ for $1\leq i\leq k$. Then we can see that the $\gamma$ with $N$ parts generates the factor $\frac{q^{(N+1)N/2}}{(q;q)_N}$ and $\delta$ generates the factor $(-1;q)_N$. We aim to show that there is a bijection between the set of partition pairs $(\gamma, \delta)$ and the set of  overpartitions $\alpha$, where the number of parts of $\alpha$ is $N$, the same as the size of its
generalized Durfee square.

First, given a partition pair $(\gamma, \delta)$, we will construct an overpartition $\alpha$. Before the construction, note the fact that $\gamma_i \geq N-i+1$, which means $\gamma_i+(i-1)\geq N$, since $\gamma$ is a partition with $N$ distinct parts.
Now we could begin to construct the overpartition $\alpha$. We overline all the parts of $\gamma$. Then from $1$ to $N$, if $\delta$ has a part $\delta_j=i-1$, add $i-1$ to $\gamma_i$ and remove the overline of $\gamma_i$.
Finally, we write the overlined parts
of $\gamma$ before the
non-overlined parts to get an overpartition  $\alpha$.
It is easy to see that
there are $N$ parts in $\alpha$ and each non-overlined part in $\alpha$ is no less than $N$.
 Then it follows that $\alpha$ is an overpartition with $N$ parts and
generalized Durfee square of size $N$.

For example, if $N=5$, $\gamma=(7,6,5,2,1)$ and $\delta=(4,3,0)$, then we have
\begin{eqnarray*}
  (7+6+5+2+1,\,4+3+0)
  &\Longrightarrow&
  (\overline{7}+\overline{6}+\overline{5}+\overline{2}+\overline{1},\,4+3+0)\\
  &\Longrightarrow&
  (7+\overline{6}+\overline{5}+\overline{2}+\overline{1},\,4+3)\\
  &\Longrightarrow&(7+\overline{6}+\overline{5}+5+\overline{1},\,4)\\
  &\Longrightarrow&(7+\overline{6}+\overline{5}+5+5)\\
  &\Longrightarrow&(\overline{6}+\overline{5}+7+5+5).
\end{eqnarray*}
Clearly, $\alpha=(\bar{6},\bar{5},7,5,5)$ have five parts and the size of generalized Durfee square is $5$.

Now we present the inverse map.
For an overpartition $\alpha=(\alpha_1, \cdots, \alpha_l, \alpha_{l+1},\cdots,\alpha_N)$, where
$\alpha_i$ is overlined for $1\leq i\leq l$ and $\alpha_j$  is non-overlined for $l+1\leq j\leq N$,
and $\alpha_N\geq N$, we construct $(\gamma,\delta)$ as follows.\\
(i) Begin with $\gamma=(\alpha_1, \cdots, \alpha_l)$ and $\delta=\varnothing$. Remove the overline of $\alpha_i$ for $1\leq i\leq l$.\\
(ii)  For $j$ from $1$ to $N-l$, choose the smallest $m\geq 0$ such that $\alpha_{l+j}-m>\gamma_{m+1}$, put $\alpha_{l+j}-m$ into $\gamma$
      as a new part, and put $m$ into $\delta$ as a new part. If $m+1>l(\gamma)$, where $l(\gamma)$ denotes the number of parts of $\gamma$, take $\gamma_{m+1}=0$.

It can be checked that $\gamma$ is a partition with $N$ distinct parts and
$\delta$ is a partition with distinct parts between $0$ and $N-1$.
And it is indeed the inverse map.

For example, if $N=5$ and $\alpha=(\overline{6},\overline{5},7,5,5)$, then we have
\begin{eqnarray*}
(\overline{6}+\overline{5}+7+5+5)
&\Longrightarrow&(6+5,\,7+5+5)\\
&\Longrightarrow&(7+6+5,\,5+5+0)\\
&\Longrightarrow&(7+6+5+2,\,5+3+0)\\
&\Longrightarrow&(7+6+5+2+1,\,4+3+0).
\end{eqnarray*}
%
%
%

This completes the proof. \qed

Now we proceed to the definition of successive Durfee squares for overpartitions. To get the successive Durfee squares, we first determine the generalized Durfee square of
the overpartition, and then the second Durfee square of the smaller partition below the generalized Durfee
square, and so on.

For example, let $\alpha=(8,\overline{7},6,6,\overline{5}, 5,5,3,\overline{1})$.
Then its successive Durfee square dissection is shown as follows.
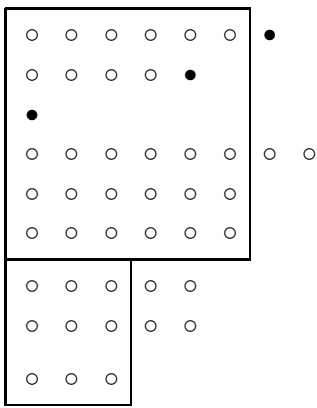
\begin{center}
 \begin{picture}(100,155)

  \put(10,10){\circle{4}}
  \put(25,10){\circle{4}}
  \put(40,10){\circle{4}}

  \put(10,30){\circle{4}}
  \put(25,30){\circle{4}}
  \put(40,30){\circle{4}}
  \put(55,30){\circle{4}}
  \put(70,30){\circle{4}}

  \put(10,45){\circle{4}}
  \put(25,45){\circle{4}}
  \put(40,45){\circle{4}}
  \put(55,45){\circle{4}}
  \put(70,45){\circle{4}}

 \put(10,65){\circle{4}}
  \put(25,65){\circle{4}}
  \put(40,65){\circle{4}}
  \put(55,65){\circle{4}}
  \put(70,65){\circle{4}}
  \put(85,65){\circle{4}}

 \put(10,80){\circle{4}}
  \put(25,80){\circle{4}}
  \put(40,80){\circle{4}}
  \put(55,80){\circle{4}}
  \put(70,80){\circle{4}}
  \put(85,80){\circle{4}}

 \put(10,95){\circle{4}}
  \put(25,95){\circle{4}}
  \put(40,95){\circle{4}}
  \put(55,95){\circle{4}}
  \put(70,95){\circle{4}}
  \put(85,95){\circle{4}}
 \put(100,95){\circle{4}}
 \put(115,95){\circle{4}}

 \put(10,110){\circle*{4.2}}

 \put(10,125){\circle{4}}
  \put(25,125){\circle{4}}
  \put(40,125){\circle{4}}
  \put(55,125){\circle{4}}
  \put(70,125){\circle*{4.2}}

\put(10,140){\circle{4}}
  \put(25,140){\circle{4}}
  \put(40,140){\circle{4}}
  \put(55,140){\circle{4}}
  \put(70,140){\circle{4}}
  \put(85,140){\circle{4}}
 \put(100,140){\circle*{4.2}}

\put(0,0){\line(0,1){150}}

\put(0,0){\line(1,0){47.5}}

\put(47.5,0){\line(0,1){55}}



\put(0,55){\line(1,0){92.5}}

\put(92.5,55){\line(0,1){95}}

\put(0,150){\line(1,0){92.5}}
 \end{picture}
\end{center}

Now we can calculate the generating function for overpartitions with at most $k-1$ successive Durfee squares
 (including the generalized Durfee square).

\begin{theo}\label{Durgf}
The generating function for overpartitions with at most $k-1$ successive Durfee squares is
\begin{equation}
 \sum_{N_1\geq N_2 \geq \ldots \geq N_k \geq 0}
\frac{q^{(N_1+1)N_1/2}q^{N_2^2+\cdots +N_{k-1}^2}(-1; q)_{N_1}}{(q;q)_{N_1-N_2}\cdots
 (q;q)_{N_{k-1}}}.
\end{equation}
\end{theo}

So the generating function of $D_{k,k}(n)$ equals the generating function for overpartitions with at most $k-1$ successive Durfee squares.

Now we consider the overpartitions with no restrictions. Recall that the  generating function of overpartition $\overline{p}(n)$ is
\begin{equation}\label{over}
  \sum_{n=0}^{\infty}\bar{p}(n)q^n= \sum_{n=0}^{\infty}\frac{(-1;q)_nq^{(n+1)n/2}}
  {(q;q)_n^2}.
\end{equation}
We can see that the Ferres diagraph  of the overpartition is composed of a generalized Durfee square with size $n$, and a partition with parts $\leq n$ under the generalized Durfee square. Obviously, the number of successive Durfee squares of the unrestricted overpartition is unrestricted.

\section{Weighted identity}

In this section, we shall construct a map between the overpartitions with at most $k-1$ successive Durfee squares and the overpartitons with the number of successive Durfee squares unrestricted. Then we can give the weight function $w_S(\pi)$ to connect $D_{k,k}(n)$ and $\overline{p}(n)$.

 To construct the map, we have to consider the gaps between parts of the generalized Durfee square, which force us to rearrange the parts in decreasing order. Meanwhile, we must make sure that the other successive Durfee squares do not change. With no loss,  we can add $N_1$ to each overlined part of the generalized Durfee square and rearrange the parts in decreasing order to get a new overpartition, where $N_1$ is the size of generalized Durfee square. And we can show that after the transformation we could minus $N_1$ from each overlined part.

Let $\lambda$ %
be an overpartition with no restriction and $\beta$
 be an overpartition with at most $k-1$ successive Durfee squares.

Now we give a surjection $\phi$: $\lambda\rightarrow \beta$. Let $k_{\lambda}$ denote the number of successive Durfee squares of $\lambda$.

(i) If  $k_{\lambda}\leq k-1$, let $\phi(\lambda)=\lambda$.

 (ii) If $k_{\lambda}\geq k$, add $N_1$ to each overlined part of the generalized Durfee square of $\lambda$ and rearrange the parts in decreasing order to get a new overpartition $\lambda'$, where $N_1$ is the size of generalized Durfee square of $\lambda$. Clearly, the number of successive Durfee squares of $\lambda'$ is $k_{\lambda}$ and the size of the first successive Durfee square of $\lambda'$ is $N_1$. Denote the partition below the $(k-1)$--th successive Durfee square of $\lambda'$ by $\lambda'_b$, and denote the partition above $\lambda'_b$ by $\lambda'_a$. Then we  minus $N_1$ from each overlined part of $\lambda'_a+(\lambda'_b)^{*}$, and get $\phi(\lambda)$ immediately.

 It is easy to see that the map $\phi$ is a surjection. Then we give the inverse map $\phi^{-1}$.

 For an overpartition $\beta$ with at most $k-1$ successive Durfee squares,  the sizes of successive Durfee squares are $N_1\geq N_2\geq \cdots \geq N_{k-1}\geq 0$.
 Let $k_{\beta}$ denote the number of successive Durfee squares of $\beta$.

 (i) If $k_{\beta}\leq k-2$,  let $\phi^{-1}(\beta)=\beta$.

 (ii) If $k_{\beta}=k-1$, we shall map $\beta$ to several overpartitions with the number of successive Durfee squares $\geq k-1$. Generally speaking, we will remove some nodes from the partition to the right of the first successive Durfee square, and place them below the bottom to construct new parts.

 step 1. Add $N_1$ to each overlined part of $\beta$ and rearrange the parts in decreasing order to get a new overpartition $\beta'$. Clearly, the first successive Durfee
  square of $\beta'$ consists of the non-overlined parts of the generalized Durfee square of $\beta$ and the overlined parts of $\beta$ with each plus $N_1$. Then each part in the first successive Durfee square of $\beta'$ is greater than or equal to $N_1$.

 Step 2. For i from $1$ to $N_{k-1}$,
 if $\beta'_{i+1}$ is overlined and $\beta'_i-\beta'_{i+1}\geq 2$, then we can subtract one  from each $\beta'_1\geq \beta'_2\geq\cdots\geq \beta'_i$, and put a part $i$ under the $(k-1)$--th successive Durfee square of $\beta'$ as a new part. Then we get a new overpartition with the number of Durfee squares $\geq k$. We can see that there are $\beta'_i-\beta'_{i+1}$ ways of the substraction.  That is, subtracting $0,\ 1,\ \ldots,$ or $\beta'_i-\beta'_{i+1}-1$ from the first $i$ parts of $\beta'$.

If $\beta'_{i+1}$ is non-overlined, and $\beta'_i-\gamma'_{i+1}\geq 1$ then we can subtract one  from each $\beta'_1\geq \beta'_2\geq\cdots\geq \beta'_i$, then we put a part $i$ under the $(k-1)$--th Durfee square of $\beta'$ as an new part. Then we can get an new overpartition with the number of durfee squares $\geq k$. We can see there are $\beta'_i-\beta'_{i+1}+1$ ways of the substraction.  That is, subtracting $0,\ 1,\ \ldots,$ or $\beta'_i-\beta'_{i+1}$ parts from the first $i$ parts of $\beta'$.

Step 3. Subtract $N_1$ from each overlined part and rearrange it.

 Then we define the weight $w(\delta)$ to be that
 \begin{equation}
 w(\delta)=\begin{cases}\prod_{i=1}^{N_{k-1}}(\beta'_i-\beta'_{i+1}+1-\varepsilon(\gamma'_{i+1})), &\text{if} \ \beta\ \text{has exactly }\ k-1\ \text{Durfee squares},\\1,&\text{if} \ \beta\ \text{has less than }\ k-1\ \text{Durfee squares}.
 \end{cases}
 \end{equation}
 where \begin{equation}\varepsilon(\beta'_{i+1})=\begin{cases}1, &\text{if }\  \beta'_{i+1}\  \text{is overlined},\\
 0, &\text{if }\ \beta'_{i+1}\  \text{is non-overlined}.
 \end{cases}
 \end{equation}

\vspace{0.5cm}
 \noindent{\bf Acknowledgments.}
This work was finished during the author Diane Y. H.
Shi's visit in RISC. We thank Peter Paule for his kindly invitation. This work  was motivated by the Krishnaswami  Alladi's talk in RISC in 2016.
This work was supported by the National Science Foundation of China (Nos.1140149, 11501089, 11501408).

\end{document}